\pdfoutput=1

\documentclass[10pt]{article}

\PassOptionsToPackage{%
	style	= alphabetic%
}{biblatex}

\PassOptionsToPackage{giveninits=true}{biblatex}

\usepackage[doi=false,isbn=false,url=false,
	hyperref=auto,
	sorting=nyt,
	maxnames=10,
	maxcitenames=3,
	backend=biber,
	texencoding=auto,
	giveninits=true,
	block=space,
]{biblatex}

\DefineBibliographyStrings{english}{%
	andothers = {et al}
} 
 
\DeclareFieldFormat{eprint:arxiv}{%
	\href{https://arxiv.org/abs/#1}{\emph{arXiv:\allowbreak #1}}}

\DeclareFieldFormat{eprint:mrnumber}{%
	\href{http://www.ams.org/mathscinet-getitem?mr=MR#1}{MR\allowbreak #1}}

\DeclareFieldFormat{eprint:jstor}{%
	\href{http://www.jstor.org/stable/#1}{JSTOR\allowbreak #1}}

\DeclareFieldFormat{eprint:customeprint}{%
	\em Available at \upshape \ttfamily \href{https://www.#1}{#1}}

\DeclareFieldFormat{eprint:online}{%
	Available \href{https://#1}{online}}

\DeclareFieldFormat{eprint:inprep}{%
	\emph{In preparation}}

\DeclareFieldFormat{eprint:manual}{%
	\emph{#1}}

\DeclareFieldFormat{eprint:onarxiv}{%
	\emph{#1}}

\DeclareFieldFormat{eprint:toappear}{%
	\mbox{\emph{To appear}}}

\DeclareFieldFormat{eprint:accepted}{%
	\emph{Accepted at {#1}; to appear}}

\DeclareSourcemap{
	\maps[datatype=bibtex]{
		\map{
			\step[fieldsource=mrnumber,		fieldtarget=eprint, final]
			\step[fieldset=eprinttype,		fieldvalue=mrnumber]
		}
		\map{
			\step[fieldsource=arxiv,		fieldtarget=eprint, final]
			\step[fieldset=eprinttype,		fieldvalue=arxiv]
		}
		\map{
			\step[fieldsource=jstor,		fieldtarget=eprint, final]
			\step[fieldset=eprinttype,		fieldvalue=jstor]
		}
		\map{
			\step[fieldsource=customeprint,	fieldtarget=eprint, final]
			\step[fieldset=eprinttype,		fieldvalue=customeprint]
		}
		\map{
			\step[fieldsource=online,		fieldtarget=eprint, final]
			\step[fieldset=eprinttype,		fieldvalue=online]
		}
		\map{
			\step[fieldsource=inprep,		fieldtarget=eprint, final]
			\step[fieldset=eprinttype,		fieldvalue=inprep]
		}
		\map{
			\step[fieldsource=manual,		fieldtarget=eprint, final]
			\step[fieldset=eprinttype,		fieldvalue=manual]
		}
		\map{
			\step[fieldsource=onarxiv,		fieldtarget=eprint, final]
			\step[fieldset=eprinttype,		fieldvalue=onarxiv]
		}
		\map{
			\step[fieldsource=toappear,		fieldtarget=eprint, final]
			\step[fieldset=eprinttype,		fieldvalue=toappear]
		}
		\map{
			\step[fieldsource=accepted,		fieldtarget=eprint, final]
			\step[fieldset=eprinttype,		fieldvalue=accepted]
		}
	}
}

\DeclareFieldFormat{doi}{%
	\href{https://doi.org/#1}{DOI}%
}

\DeclareFieldFormat{url}{%
	\href{https://#1}{\texttt{#1}}%
}

\DeclareFieldFormat{comments}{%
	\emph{#1}%
} 
 
\DeclareFieldFormat
	[article,inbook,incollection,inproceedings,patent,thesis,unpublished]
	{title}{{#1\isdot}}

\DeclareFieldFormat
	[article,book,inbook,incollection,inproceedings,patent,thesis,unpublished, online]
	{date}{\mkbibparens{#1}}
\DeclareFieldFormat
	[article,book,inbook,incollection,inproceedings,patent,thesis,unpublished]
	{volume}{\mkbibbold{#1}}
\DeclareFieldFormat
	{pages}{\mkbibparens{#1}}

\DeclareFieldFormat{edition}{%
	\ifinteger{#1}
	{\mkbibordedition{#1}~\bibstring{edition}\addcomma}
	{#1~\bibstring{edition}\addcomma}}

\renewbibmacro*{volume+number+eid}{%
	\printfield{volume}%
\setunit*{\adddot}%
	\printfield{number}%
\setunit{\space}%
	\printfield{eid}%
}

\DeclareBibliographyDriver{article}{%
	\usebibmacro{bibindex}%
	\usebibmacro{begentry}%
	\usebibmacro{author}%
\setunit{\addspace}%
	\usebibmacro{date}%
\newunit\newblock
	\usebibmacro{title}%
\newunit\newblock
	\printfield{journaltitle}\addperiod%
\setunit*{\addspace}%
	\usebibmacro{volume+number+eid}%
\setunit*{\addspace}\newblock
	\printfield{pages}
\setunit*{\addspace}
	\printfield{comments}%
	\usebibmacro{eprint}%
\setunit{\addnbspace}%
	\printfield{doi}%
	\usebibmacro{finentry}%
}

\DeclareBibliographyDriver{online}{%
	\usebibmacro{bibindex}%
	\usebibmacro{begentry}%
	\usebibmacro{author}%
\setunit{\addspace}%
	\usebibmacro{date}%
\newunit\newblock
	\usebibmacro{title}%
\newunit\newblock
	\usebibmacro{eprint}%
	\usebibmacro{finentry}%
}

\DeclareBibliographyDriver{book}{%
	\usebibmacro{bibindex}%
	\usebibmacro{begentry}%
	\usebibmacro{author}%
\setunit{\addspace}\newblock
	\usebibmacro{date}%
\newunit\newblock
	\usebibmacro{title}%
\setunit*{\addspace}\newblock
	\printfield{edition}%
\newunit\newblock
	\printlist{publisher}
\setunit*{\addspace}%
	\printfield{volume}%
\setunit*{\addspace}\newblock%
	\printfield{comments}%
	\usebibmacro{eprint}%
\setunit{\addnbspace}
	\printfield{doi}
	\usebibmacro{finentry}%
}

\DeclareBibliographyDriver{thesis}{%
	\usebibmacro{bibindex}%
	\usebibmacro{begentry}%
	\usebibmacro{author}%
\setunit{\addspace}\newblock
	\usebibmacro{date}%
\newunit\newblock
	\usebibmacro{title}.%
\setunit*{\addspace}\newblock
	\space Thesis, \printlist{institution}
	\usebibmacro{eprint}%
	\printfield{comments}%
\setunit*{\addnbspace}
	\printfield{doi}%
	\usebibmacro{finentry}%
}

\DeclareBibliographyDriver{inproceedings}{%
	\usebibmacro{bibindex}%
	\usebibmacro{begentry}%
	\usebibmacro{author}%
\setunit{\addspace}%
	\usebibmacro{date}%
\newunit\newblock
	\usebibmacro{title}%
\newunit\newblock
	\printfield{booktitle}\addcomma%
\setunit*{\addspace}%
	\printfield{series}\addcomma%
\setunit*{\addspace}\newblock
	\usebibmacro{editor}\addcomma
\setunit*{\addspace}\newblock
	\printlist{publisher}
\setunit*{\addspace}%
	\printfield{volume}%
\setunit*{\addspace}%
	\printfield{pages}
\setunit*{\addspace}
	\usebibmacro{eprint}%
	\printfield{comments}%
\setunit*{\addnbspace}
	\printfield{doi}
	\usebibmacro{finentry}%
}

\DeclareBibliographyDriver{inbook}{%
	\usebibmacro{bibindex}%
	\usebibmacro{begentry}%
	\usebibmacro{author}%
\setunit{\addspace}%
	\usebibmacro{date}%
\newunit\newblock
	\usebibmacro{title}%
\newunit\newblock
	\printfield{booktitle}\addcomma%
\setunit*{\addspace}%
	\printfield{series}\addcomma%
\setunit*{\addspace}\newblock
	\usebibmacro{editor}\addcomma
\setunit*{\addspace}\newblock
	\printlist{publisher}
\setunit*{\addspace}%
	\printfield{volume}%
\setunit*{\addspace}%
	\printfield{pages}
\setunit*{\addspace}
	\usebibmacro{eprint}%
	\printfield{comments}%
\setunit*{\addspace}
	\printfield{doi}
	\usebibmacro{finentry}%
}

\DeclareBibliographyDriver{incollection}{%
	\usebibmacro{bibindex}%
	\usebibmacro{begentry}%
	\usebibmacro{author}%
\setunit{\addspace}%
	\usebibmacro{date}%
\newunit\newblock
	\usebibmacro{title}%
\newunit\newblock
	\printfield{booktitle}\addcomma%
\setunit*{\addspace}%
	\printfield{series}\addcomma%
\setunit*{\addspace}\newblock
	\printlist{publisher}
\setunit*{\addspace}%
	\printfield{volume}%
\setunit*{\addspace}%
	\printfield{pages}
\setunit*{\addspace}
	\usebibmacro{eprint}%
	\printfield{comments}%
\setunit{\addnbspace}
	\printfield{doi}
	\usebibmacro{finentry}%
}

\AtEveryBibitem{%
	\clearfield{day}%
	\clearfield{month}%
} 
 
\DeclareNameAlias{sortname}{given-family}
\DeclareNameAlias{default}{given-family}

\usepackage{geometry}
\newcommand{\MCH}{\ifmmode \mathsf{MCH} \else \textsf{MCH}\xspace \fi}
\newcommand{\QMCH}{\ifmmode \mathsf{QMCH} \else \textsf{QMCH}\xspace \fi}
\newcommand{\TV}{\ifmmode \mathsf{TV} \else \textsf{TV}\xspace \fi}

\usepackage[utf8]{inputenc}

\usepackage{chngcntr}

\usepackage[UKenglish]{babel}
\usepackage{amsmath}
\usepackage{amsthm}
\usepackage{amssymb}
\usepackage{mathrsfs}
\usepackage{bm}
\usepackage{mathtools}
\newcommand{\cq}{\coloneqq}

\usepackage{titlesec}
\usepackage{xifthen}
\usepackage{xspace}

\usepackage{graphicx}
\usepackage[labelfont = {bf, sf}, labelsep = period]{caption}
\usepackage{subcaption}

\usepackage{titletoc}

\titleformat{\section}
	{\sffamily \Large \bfseries \boldmath}{\thesection}{1em}{}
\titleformat{\subsection}
	{\sffamily \large \bfseries \boldmath}{\thesubsection}{0.5em}{}

\def\IfAmpersandUseAlign#1#2&#3\EndIfAmpersandUseAlign%
{%
	\if\relax\detokenize{#3}\relax
	\begin{equation*}%
	#1%
	\end{equation*}%
	\else
	\begin{align*}%
	#1%
	\end{align*}%
	\fi
}
\def\[#1\]%
{%
	\IfAmpersandUseAlign{#1}#1&\EndIfAmpersandUseAlign
}

\newcommand{\off}{\textsc{off}\xspace}
\newcommand{\on}{\textsc{on}\xspace}

\usepackage{manyfoot}
\SetFootnoteHook{\hspace*{-1.8em}}
\DeclareNewFootnote{bl}[gobble]
\newcommand{\blfootnote}[1]{\footnotebl{#1}}

\newcommand{\nt}{\addtocounter{equation}{1}\tag{\theequation}}

\newcommand{\eps}{\varepsilon}
\newcommand{\bcdot}{\ensuremath{\bm{\cdot}}}

\newcommand{\Quad}[1]{\quad \text{#1} \quad}

\newcommand{\mix}{\textnormal{mix}}

\newcommand{\tmix}{t_\mix}

\DeclareMathOperator{\Unif}{Unif}
\newcommand{\iid}{\textup{iid}}

\DeclareMathOperator{\Bin}{Bin}

\DeclareMathOperator{\Pois}{Pois}
\DeclareMathOperator{\Bern}{Bern}

\usepackage{enumitem}

\setlist[itemize,enumerate]{%
	itemsep	= 0pt,%
	topsep	= \smallskipamount,%
	label	= \ensuremath{\bcdot}%
}
\setlist[description]{%
	topsep = \smallskipamount,		
	itemsep = \smallskipamount,		
	font = {\mdseries\itshape},		
	leftmargin = \parindent
}

\usepackage{xifthen}

\newcommand{\maxt}[1]{ \textstyle \max_{#1} \displaystyle }

\newcommand{\sumt}[2][]{
	\ifthenelse{\isempty{#1}}
	{\textstyle \sum_{#2}      \displaystyle}
	{\textstyle \sum_{#2}^{#1} \displaystyle}
}
\newcommand{\sumd}[2][]{
	\ifthenelse{\isempty{#1}}
	{\sum_{#2}}
	{\sum_{#2}^{#1}}
}

\newcommand{\intt}[2][]{
	\ifthenelse{\isempty{#1}}
	{\textstyle \int_{#2}      \displaystyle}
	{\textstyle \int_{#2}^{#1} \displaystyle}
}
\newcommand{\intd}[2][]{
	\ifthenelse{\isempty{#1}}
	{\int_{#2}}
	{\int_{#2}^{#1}}
}

\newcommand{\prodt}[2][]{
	\ifthenelse{\isempty{#1}}
	{\textstyle \prod_{#2}      \displaystyle}
	{\textstyle \prod_{#2}^{#1} \displaystyle}
}
\newcommand{\prodd}[2][]{
	\ifthenelse{\isempty{#1}}
	{\prod_{#2}}
	{\prod_{#2}^{#1}}
}

\newcommand{\pr}[2][]{
	\ifthenelse{\equal{}{#1}}
	{\mathbb{P}[#2]}
	{\mathbb{P}_{#1}[#2]}
}
\newcommand{\prb}[2][]{
	\ifthenelse{\equal{}{#1}}
	{\mathbb{P}\bigl[ #2 \bigr]}
	{\mathbb{P}_{#1}\bigl[ #2 \bigr]}
}
\newcommand{\prB}[2][]{
	\ifthenelse{\equal{}{#1}}
	{\mathbb{P}\Bigl[ #2 \Bigr]}
	{\mathbb{P}_{#1}\Bigl[ #2 \Bigr]}
}
\newcommand{\prbb}[2][]{
	\ifthenelse{\equal{}{#1}}
	{\mathbb{P}\biggl[ #2 \biggr]}
	{\mathbb{P}_{#1}\biggl[ #2 \biggr]}
}
\newcommand{\prBB}[2][]{
	\ifthenelse{\equal{}{#1}}
	{\mathbb{P}\Biggl[ #2 \Biggr]}
	{\mathbb{P}_{#1}\Biggl[ #2 \Biggr]}
}
\newcommand{\prs}[2][]{
	\ifthenelse{\equal{}{#1}}
	{\mathbb{P}\left[ #2 \right]}
	{\mathbb{P}_{#1}\left[ #2 \right]}
}

\newcommand{\ex}[2][]{
	\ifthenelse{\equal{}{#1}}
	{\mathbb{E}[#2]}
	{\mathbb{E}_{#1}[#2]}
}
\newcommand{\exb}[2][]{
	\ifthenelse{\equal{}{#1}}
	{\mathbb{E}\bigl[ #2 \bigr]}
	{\mathbb{E}_{#1}\bigr[ #2 \bigr]}
}
\newcommand{\exB}[2][]{
	\ifthenelse{\equal{}{#1}}
	{\mathbb{E}\Bigl[ #2 \Bigr]}
	{\mathbb{E}_{#1}\Bigl[ #2 \Bigr]}
}
\newcommand{\exbb}[2][]{
	\ifthenelse{\equal{}{#1}}
	{\mathbb{E}\biggl[ #2 \biggr]}
	{\mathbb{E}_{#1}\biggl[ #2 \biggr]}
}
\newcommand{\exBB}[2][]{
	\ifthenelse{\equal{}{#1}}
	{\mathbb{E}\Biggl[ #2 \Biggr]}
	{\mathbb{E}_{#1}\Biggl[ #2 \Biggr]}
}

\newcommand{\Varb}[2][]{
	\ifthenelse{\equal{}{#1}}
	{\mathbb{V}\textnormal{ar} \bigl[#2\bigr]}
	{\mathbb{V}\textnormal{ar}_{#1} \bigl[#2\bigr]}
}
\newcommand{\VAR}[2][]{
	\ifthenelse{\equal{}{#1}}
	{\textnormal{Var}[#2]}
	{\textnormal{Var}_{#1}[#2]}
}

\newcommand{\tv} [1]{\| #1 \|_\TV}
\newcommand{\tvb}[1]{\bigl\| #1 \bigr\|_\TV}

\newcommand{\abs}  [1]{\lvert #1 \rvert}
\newcommand{\absb} [1]{\big\lvert #1 \bigr\rvert}

\newcommand{\norm}  [1]{\lVert #1 \rVert}

\newcommand{\rbr} [1]{ ( #1 ) }
\newcommand{\rbb} [1]{\bigl( #1 \bigr)}

\newcommand{\sbr} [1]{ [ #1 ] }
\newcommand{\sbb} [1]{\bigl[ #1 \bigr]}

\newcommand{\bra} [1]{ \{ #1 \} }
\newcommand{\brb} [1]{\bigl\{ #1 \bigr\}}

\newcommand{\floor}[1]{\lfloor #1 \rfloor}

\newcommand{\one}  [1]{\bm1\{ #1 \}}

\let\originalexp\exp
\let\exp\relax
\DeclareRobustCommand{\exp} [1]{\originalexp(#1)}

\newcommand{\oh}  [1]{o( #1 )}

\newcommand{\mbn}{\ensuremath{\mathbb{N}}}

\newcommand{\mbr}{\ensuremath{\mathbb{R}}}

\newcommand{\mbz}{\ensuremath{\mathbb{Z}}}

\newcommand{\mca}{\ensuremath{\mathcal{A}}}

\newcommand{\mck}{\ensuremath{\mathcal{K}}}

\usepackage[dvipsnames,hyperref]{xcolor}


\newcounter{parentnumber}

\makeatletter
\newenvironment{subtheorem-num}[1]{%
	\def\subtheoremcounter{#1}%
	\refstepcounter{#1}%
	\protected@edef\theparentnumber{\csname the#1\endcsname}%
	\setcounter{parentnumber}{\value{#1}}%
	\setcounter{#1}{0}%
	\expandafter\def\csname the#1\endcsname{\theparentnumber.\arabic{#1}}%
	\expandafter\def\csname theH#1\endcsname{thm.\theparentnumber.\arabic{#1}}%
	\unskip\ignorespaces
}{%
	\setcounter{\subtheoremcounter}{\value{parentnumber}}%
	\ignorespacesafterend
}
\makeatother

\usepackage[
	colorlinks,
	pdfauthor = {Sam\ Olesker-Taylor},
	pdfusetitle
]{hyperref}
\hypersetup{
	colorlinks,
	linkcolor	= {blue},
	urlcolor	= {blue},
	citecolor	= {ForestGreen}
}

\usepackage[capitalise]{cleveref}

\newenvironment{Proof}[1][\proofname]{%
	\proof[\upshape\bfseries\sffamily\boldmath#1]
}{\endproof}

\newtheoremstyle{sfsl}
{1\baselineskip}		
{1\baselineskip}		
{\slshape}				
{}						
{\bfseries\sffamily}	
{.}						
{0.5em}					
{\thmname{#1}\thmnumber{ #2}\thmnote{ \textnormal{\sffamily(#3)}}}

\newtheoremstyle{sfup}
{1\baselineskip}		
{1\baselineskip}		
{\upshape}				
{}						
{\bfseries\sffamily}	
{.}						
{0.5em}					
{\thmname{#1}\thmnumber{ #2}\thmnote{ \textnormal{\sffamily(#3)}}}

\newtheoremstyle{bvsl}
{1\baselineskip}		
{1\baselineskip}		
{\slshape}				
{}						
{\bfseries}				
{.}						
{0.5em}					
{\thmname{#1}\thmnumber{ #2}\thmnote{ \textnormal{(#3)}}}

\newtheoremstyle{bvup}
{1\baselineskip}		
{1\baselineskip}		
{\upshape}				
{}						
{\bfseries}				
{.}						
{0.5em}					
{\thmname{#1}\thmnumber{ #2}\thmnote{ \textnormal{(#3)}}}

\theoremstyle{sfsl}

\newtheorem*{thm*}{Theorem}

\crefname{thm}{Theorem}{Theorems}

\newtheorem*{introthm*}{Theorem}
\newtheorem{introthm}{Theorem}

\crefname{introthm}{Theorem}{Theorems}

\newtheorem{introcor}{Corollary}

\crefname{introcor}{Corollary}{Corollaries}

\newtheorem{introprop}{Proposition}

\crefname{introprop}{Proposition}{Propositions}

\newtheorem*{prop*}    {Proposition}

\crefname{prop}{Proposition}{Propositions}

\newtheorem*{defn*}    {Definition}

\crefname{defn}{Definition}{Definitions}

\newtheorem*{introdefn*}{Definition}
\newtheorem{introdefn}{Definition}

\crefname{introdefn}{Definition}{Definitions}

\newtheorem{innercustomgeneric}{\customgenericname}
\providecommand{\customgenericname}{}
\newcommand{\newcustomtheorem}[2]{%
	\newenvironment{#1}[1]
	{%
		\renewcommand\customgenericname{#2}%
		\renewcommand\theinnercustomgeneric{##1}%
		\innercustomgeneric
	}
	{\endinnercustomgeneric}
}

\newcustomtheorem{customthm}{Theorem}
\newcustomtheorem{customprop}{Proposition}
\newcustomtheorem{customlemma}{Lemma}
\newcustomtheorem{customconj}{Conjecture}
\newcustomtheorem{customquestion}{Question}
\newcustomtheorem{customopenproblem}{Open Problem}

\newtheorem*{conj*}   {Conjecture}

\crefname{conj}{Conjecture}{Conjectures}

\newenvironment{conj-ind*}
	{\begin{quote}\textsf{\textbf{Conjecture.}}\slshape}
	{\end{quote}}
\newenvironment{conj-ind}
	{\begin{quote}\vspace{-\glueexpr\baselineskip+\topsep}\begin{customconj}}
	{\end{customconj}\end{quote}}

\newenvironment{question-ind*}
	{\begin{quote}\textsf{\textbf{Question.}}\slshape}
	{\end{quote}}
\newenvironment{question-ind}
	{\begin{quote}\vspace{-\glueexpr\baselineskip+\topsep}\begin{customquestion}}
	{\end{customquestion}\end{quote}}

\newenvironment{openproblem-ind*}
	{\begin{quote}\textsf{\textbf{Open Problem.}}\slshape}
	{\end{quote}}
\newenvironment{openproblem-ind}
	{\begin{quote}\vspace{-\glueexpr\baselineskip+\topsep}\begin{customopenproblem}}
	{\end{customopenproblem}\end{quote}}

\newtheorem*{hypothesis*}{Hypothesis}

\newtheorem*{hyp*}{Hypotheses}

\crefname{hyp}{Hypotheses}{Hypotheses}

\theoremstyle{sfup}

\crefname{rmk} {Remark}{Remarks}
\crefname{rmkT}{Remark}{Remarks}

\newtheorem*{rmk*} {Remark}
\newtheorem*{rmkTT}{Remark}
\newenvironment{rmkt*}
	{\pushQED{\qed}\rmkTT}
	{\popQED\endrmkTT}

\crefname{rmks} {Remarks}{Remarks}
\crefname{rmksT}{Remarks}{Remarks}

\newtheorem*{rmks*} {Remarks}
\newtheorem*{rmksTT}{Remarks}
\newenvironment{rmkst*}
	{\pushQED{\qed}\rmksTT}
	{\popQED\endrmksTT}

\newtheorem{intrormkT}{Remark}
	\newenvironment{intrormkt}
	{\pushQED{\qed}\intrormkT}
	{\popQED\endintrormkT}

\crefname{intrormk} {Remark}{Remarks}
\crefname{intrormkT}{Remark}{Remarks}

\newtheorem*{intrormk*} {Remark}
\newtheorem*{intrormkTT}{Remark}
\newenvironment{intrormkt*}
	{\pushQED{\qed}\intrormkTT}
	{\popQED\endintrormkTT}

\newcommand{\BM}{}

\title{\sffamily%
	Multicoloured Hardcore Model: Fast Mixing and Queueing}
\author{\sffamily%
	Sam Olesker-Taylor}
\date{\sffamily%
	\today}
\date{}

\begin{document}

\maketitle

\renewcommand{\abstractname}{\sffamily Abstract}

\begin{abstract}
\noindent
In the hardcore model, certain vertices in a graph are \textit{active}: these must form an \textit{independent set}. We extend this to a \textit{multicoloured} version: instead of simply being active or not, the active vertices are assigned a colour; active vertices \emph{of the same colour} must not be adjacent.

This models a scenario in which two neighbouring resources may \emph{interfere} when active---eg, short-range radio communication. However, there are multiple \emph{channels} (colours) available; they only interfere if both use the \emph{same} channel. Other applications include fibreoptic routing.

We analyse Glauber-type dynamics. A vertex updates its status at the incidents of a Poisson process, at which a biased coin is tossed and a uniform colour is proposed: the vertex is assigned that colour if the coin comes up heads and the colour is available; otherwise, it~is~deactivated.

We find conditions for \textit{fast mixing} of these dynamics. We also use them to model a queueing system: vertices only serve their customers whilst active. The mixing estimates are applied to establish positive recurrence of the queue lengths, and bound their expectation in equilibrium.
\end{abstract}

\small
\begin{quote}
\begin{description}
	\item [Keywords:]
	mixing time,
	hardcore model,
	proper colourings,
	queueing theory,
	scheduling
	
	\item [MSC 2020 subject classifications:]
	60C05;
	60J20,
	60J27;
	60K25,
	60K30;
	68M20
\end{description}
\end{quote}
\normalsize







\addtocontents{toc}{\vspace{-\medskipamount}}
\sffamily
\setcounter{tocdepth}{1}
\tableofcontents
\normalfont

\blfootnote{\sffamily%
	Department of Statistics, University of Warwick, UK
\hfill%
	\href{mailto:oleskertaylor.sam@gmail.com}{oleskertaylor.sam@gmail.com}%
}


\vspace{-\medskipamount}

\renewcommand{\thesubsection}{\Alph{subsection}}
\section{Introduction and Main Results}

We extend the well-studied hardcore model,
used for sampling independent sets,
to a multicoloured version.
More precisely, given a graph $G = (V, E)$, our objective is to colour a subset $U \subseteq V$ of the vertices such that if $u,u' \in U$ satisfy $\bra{u,u'} \in E$, then $u$ and $u'$ are painted with different~colours.
If there is only one colour, then this condition requires that there is no pair of mutually adjacent vertices.
This is the definition of an \textit{independent set}, so we recover the usual \textit{hardcore model}.

We allow an arbitrary number $K \in \mbn$ of colours.
If we required \emph{all} vertices to be selected---ie, $U = V$---then the condition is that no edge in the graph is \textit{monochromatic}: the endpoints must receive different colours.
We thus recover the \textit{proper colouring} model, which is also well-studied.
Our model models these two, sampling a properly coloured subset of vertices, or~subgraph.

The motivation for this model comes from a desire for a decentralised (and randomised) algorithm for resource sharing.
Two examples of this
are short-range radio communication, where nearby agents on the same frequency interfere, and routing algorithms in fibreoptic networks.
Both $K$ and $G$ are given parameters, depending on the particular engineering set-up.

\medskip

A popular method for sampling proper colourings or independent sets is via \textit{Glauber dynamics}.
Our main result is on the \textit{mixing time} of Glauber-type dynamics for the multicoloured hardcore model, defined precisely below.
We then use the system to model a \textit{queueing system}.

\begin{itemize}
	\item 
	Customers (eg, data packets) arrive to vertices at some (vertex-dependent) rate.
	
	\item 
	Coloured vertices are \textit{active}: they serve their customers at some (vertex-dependent) rate.
	
	\item 
	Uncoloured vertices are \textit{inactive}: they do not serve, but their queue can still grow.
\end{itemize}
We apply the mixing-time result to control queue lengths of this system,
under certain conditions.

\begin{subtheorem-num}{introdefn}
\label{sub:intro:mix}

\subsection{Glauber Dynamics}

Let $G = (V, E)$ be a graph and $K \in \mbn$.
Let $n \cq \abs V$;
write
	$[K]_0 \cq \bra{0, 1, ..., K}$.
The state space $\Omega$ of the system is a subset of configurations
\(
	[K]_0^V
=
	\bra{ \BM \omega = (\omega_v)_{v \in V} \mid \omega_v \in [K]_0 \: \forall \, v \in V }.
\)

\begin{introdefn}[State Space]
Let
\(
	\Omega
\cq
	\bra{ \BM \omega \in [K]_0^V \mid \text{$\BM \omega$ is proper} },
\)
where $\BM \omega \in [K]_0^V$ is \textit{proper} if
\[
	\omega_u \ne \omega_v
\Quad{whenever}
	\bra{u,v} \in E
\Quad{and}
	\omega_u + \omega_v > 0.
\]
\end{introdefn}

For a configuration to be proper, the colour of one vertex must be different to that of all its neighbours, except that colour $0$ is exempt from this condition.
We think of colour $0$ as inactive.
Then, a configuration is proper if the subgraph induced by its active vertices is properly coloured.


\begin{introdefn}[Glauber-Type Dynamics]
Let $\BM \lambda = (\lambda_v) \in (0,\infty)^V$ and $\BM p = (p_v) \in [0,1]^V$.
We analyse the following continuous-time Markov chain on $\Omega$, which we denote $\MCH_\Omega(\BM \lambda, \BM p)$.

\begin{itemize}
	\item 
	Select vertex $v \in V$ to update at rate $\lambda_v$, simultaneously over all vertices.
	
	\item 
	When vertex $v \in V$ is selected, toss a coin $C \sim \Bern(p_v)$ and draw a colour $k \sim \Unif([K])$.
	\begin{itemize}
		\item 
		If $C = 1$ and
			no neighbour of $v$ has colour $k$
			(ie, $k$ is \textit{available}),
		then set $v$ to colour $k$.
		
		\item 
		Otherwise%
			---ie, if $C = 0$ or the proposed colour is unavailable---%
		set $v$ to colour $0$.
	\end{itemize}
\end{itemize}
Denote the equilibrium distribution by $\pi$.
The \textit{equilibrium service rates} $\BM s = (s_v)_{v \in V}$ are given by
\[
	s_v
\cq
	\sumt{\omega \in \Omega : \omega_v \ne 0}
	\pi(\omega)
\Quad{for}
	v \in V.
\]
\end{introdefn}

The usual Glauber dynamics for proper colourings proposes a colour chosen uniformly \emph{amongst available colours}.
However, this requires whoever is making the colour choice to know which colours are available for that vertex.
This is unreasonable in the context of routing algorithms in fibreoptic networks, for example.
It is often much faster to check if a single proposed colour is available than to find out which colours are available.
We discuss this further in the motivation section.

Our main theorem regards \emph{fast mixing} of these dynamics.
First, we define mixing times precisely.

\begin{introdefn}[Mixing Times]
The \textit{total-variation} distance
between distributions $\mu$ and $\pi$ on $\Omega$~is
\[
	\tv{ \mu - \pi }
\cq
	\maxt{A \subseteq \Omega} \,
	\abs{ \mu(A) - \pi(A) }
=
	\tfrac12
	\sumt{\BM \omega \in \Omega}
	\abs{ \mu(\BM \omega) - \pi(\BM \omega) }.
\]
The \textit{mixing time} of a Markov chain $X = (X^t)_{t\ge0}$ on $\Omega$ with invariant distribution $\pi$ is
\[
	\tmix(\eps)
\cq
	\inf\bra{ t \ge 0 \mid \maxt{x \in \Omega} \tv{ \pr[x]{X^t \in \cdot} - \pi } \le \eps }
\Quad{for}
	\eps \in (0,1).
\]
\end{introdefn}

\end{subtheorem-num}


\begin{introthm}[Fast Mixing]
\label{thm:mix}
Suppose that there exists $\beta > 0$ such that
\[
	\tfrac1K
	\sumt{u \in V : \bra{u,v} \in E}
	p_u \lambda_u / \lambda_v
\le
	1 - \beta
\Quad{for all}
	v \in V.
\]
Let $\lambda_{\min} \cq \min_{v \in V} \lambda_v$.
If $\BM X, \BM Y \sim \MCH_\Omega(\BM \lambda, \BM p)$, then
\[
	\max_{\BM x, \BM y \in \Omega} \,
	\tvb{ \pr[\BM x]{ \BM X^t \in \cdot } - \pr[\BM y]{ \BM Y^t \in \cdot } }
\le
	\min\brb{ 2n e^{- \beta \lambda_{\min} t}, \: 1 }.
\]
In particular,
\[
	\tmix(\eps)
\le
	(\beta \lambda_{\min})^{-1}
	\log(2n/\eps)
\Quad{for all}
	\eps \in (0,1).
\]
\end{introthm}

\begin{intrormkt}[Fast-Mixing Condition]
The condition in \cref{thm:mix} arises from requiring the Wasserstein distance between $\BM X$ and $\BM Y$ to contract in a single step, uniformly.
Distance is measured vertex-wise: $d(\BM x, \BM y) \cq \sum_{v \in V} \one{x_v \ne y_v}$ for $\BM x, \BM y \in \Omega$.
Namely, we prove that
if configurations $\BM x$ and $\BM y$ differ only in that vertex $v$ is active in one but not the other, then
\[
	\left.
		\tfrac d{dt}
		\ex[(\BM x, \BM y)]{ d(\BM X^t, \BM Y^t) }
	\right|_{t=0}
\le
	\lambda_v
	\rbb{ \tfrac1K \sumt{u \in V : \bra{u,v} \in E} p_u \lambda_u / \lambda_v - 1 }
\]
under some coupling. The condition in \cref{thm:mix} ensures this is negative, uniformly in $v$.
A standard application of \textit{path coupling}~\cite{BD:path-coupling} extends this uniform contraction to all $\BM x, \BM y \in \Omega$.
\end{intrormkt}

The graph $G$ and number $K$ of colours (eg, channels) are given by the application.
In contrast, the parameters $(\lambda_v, p_v)_{v \in V}$ may be chosen by the operator.
There are good heuristics for taking
\[
	\lambda_v \propto d_v
\Quad{and}
	p_v \propto (K/d_v) \wedge 1.
\]
In short, high-degree nodes have more impact on their neighbours, and hence should be updated faster: so, take $\lambda_v \propto d_v$.
Further, if $v$ is active with probability $s_v$, then it removes a total of $s_v d_v$ colour choices in expectation (from its neighbours).
There are $K$ colours, so vertices shouldn't remove more than $K$ in expectation: hence, $s_v d_v \propto K$; so, take $p_v \propto (K/d_v) \wedge 1$.

We work in continuous time, so scaling all the rates $\lambda$ inversely scales the mixing time by the same amount.
We choose the normalisation $\sum_{v \in V} \lambda_v = n$;
so, vertices update at rate $1$ on average.

\begin{introcor}[Heuristic-Driven Choice]
\label{cor:mix}
Suppose that $\lambda_v = d_v / \bar d$ and $p_v \le \tfrac23 K / d_v$ for all $v \in V$, where $\bar d \cq \tfrac1n \sum_{v \in V} d_v$ is the average degree.
Let $\delta \cq \min_{v \in V} d_v$.
Let $X, Y \sim \MCH_\Omega(\BM \lambda, \BM p)$.~%
Then,
\[
	\max_{x,y \in \Omega}
	\tvb{ \pr[x]{ X^t \in \cdot } - \pr[y]{ Y^t \in \cdot } }
\le
	\min\brb{ 2 n e^{- (\delta/\bar d) t / 3}, \: 1 }.
\]
In particular,
\[
	\tmix(\eps)
\le
	3 (\bar d / \delta) \log(2n/\eps)
\Quad{for all}
	\eps \in (0,1).
\]
\end{introcor}

It is standard, or, at least, very common, in the hardcore-model ($K = 1$) literature to require $p_v = p < 1/\Delta$ for all $v$,
where $\Delta \cq \max_{v \in V} d_v$ is the maximum degree;
see, eg, \cite{BST:hardcore,LYZ:hardcore,GGSVY:hardcore} or \cite[Theorem~5.9]{LPW:markov-mixing}.
We take more care, requiring only $p_v < K/d_v$ for each $v$;
\cite{JLNRW:parallel} have a similar improvement, but restricted to the usual hardcore model ($K = 1$).

A consequence of requiring $p_v = p < 1/\Delta$ is that the mixing time is often proportional to~$\Delta$.
Ours is proportional to $\bar d / d_{\min}$,
which is often significantly smaller.

The bound $p < 1/\Delta$ is natural, up to a factor $e$.
Indeed, for the (usual, single-colour) hardcore model, it has been known since \textcite{K:stochastic-comms} that the infinite $\Delta$-regular tree has a critical threshold at $p_c(\Delta) \approx e/\Delta$, for large $\Delta$:
	the corresponding Gibbs distribution is unique if and only if $p < p_c(\Delta)$.
When $p < p_c(\Delta)$, known as the \textit{uniqueness regime} to physicists, the `influence' of one vertex on another decays exponentially in their relative distance. On the other hand, long-rage correlations persist when $p > p_c(\Delta)$.
More discussion on this can be found in \cite[\S 1.2]{ALOg:spectral-indep}.

Based on this, it appears that we should be able to only require $p_v \le (1 - \eta) e K / d_v$ and still obtain fast mixing.
This would be a natural extension of the critical threshold: $p_c(\Delta, K) \cq e K / \Delta$.
We demonstrate this via some simulations in the final section of the paper.



\medskip

We also investigate the proportion of time that vertices are active in equilibrium.

\begin{introprop}[Equilibrium Service Rates]
\label{prop:rates}
Suppose that $p_v \le \tfrac13 K / \tilde d_v$ for all $v \in V$,
where $\tilde d_v \cq \max\bra{ d_u \mid u \sim v \text{ or } u = v }$ is the maximal degree in the neighbourhood of $v \in V$.
Then,
\[
	\tfrac13
	p_v
\le
	s_v
\le
	p_v
\Quad{for all}
	v \in V.
\]
\end{introprop}


Our proof is quite flexible, not requiring these specific bounds on $p_v$.
We discuss how to generalise it, and tighten the bound, after its proof.
Again, we expect that we only need $p_v \le (1 - \eta) e K / d_v$.


\subsection{Queueing Network}

Our second results concerns queue length in a queueing network.
The proof relies on fast mixing.

\begin{introdefn}[Queueing Network]
	%
Let $\BM \lambda, \BM \nu, \BM \mu \in (0, \infty)^V$ and $\BM p \in [0,1]^V$.
Let $X \sim \MCH_\Omega(\BM \lambda, \BM p)$.
The state space of the \textit{queueing network} is $\mbn^V$.
For $q \in \mbn^V$ and $v \in V$, let
\[
	q^{v,\pm}_u
\cq
	q_u \pm \one{u = v}
\Quad{for}
	u \in V;
\]
that is, $q^{v,\pm}$ adds/removes one from the $v$-th queue.
The transition rates given $X = x$ are
\[
	q
\to
\begin{cases}
	q^{v,+}
		&\text{at rate}\quad
	\nu_v
\\
	q^{v,-}
		&\text{at rate}\quad
	\mu_v \one{x_v \ne 0}
\end{cases}
\Quad{for each}
	v \in V;
\]
that is, the $v$-th queue has arrivals at rate $\nu_v$ and, \emph{provided} $v$ is active, services at rate $\mu_v$.
We denote the law of this queueing network by $\QMCH_\Omega(\BM \lambda, \BM p; \BM \nu, \BM \mu)$.
We take $\mu_v \cq 1$ for all $v \in V$ below.
\end{introdefn}


We show that the queues are jointly \textit{positive recurrent}%
	---ie, the expected time until all queues are simultaneously empty is finite---%
under the fast-mixing conditions of \cref{thm:mix} and the assumption that the arrival rate $\nu_v$ is smaller than the equilibrium service rate $s_v$ for all $v \in V$.


\begin{introthm}[Stable Queues]
\label{thm:q}
Suppose that there exists $\beta > 0$ such that
\[
	\tfrac1K
	\sumt{u \in V : \bra{u,v} \in E}
	p_u \lambda_u / \lambda_v
\le
	1 - \beta
\Quad{for all}
	v \in V.
\]
Suppose also that $\nu_v < s_v$ for all $v \in V$.
If $Q \sim \QMCH_\Omega(\BM \lambda, \BM p; \BM \nu, \BM 1)$, then $Q$ is positive recurrent:
\[
	\tau
\cq
	\inf\bra{ t \ge 0 \mid Q^t = 0 }
\Quad{satisfies}
	\ex[q]{\tau} < \infty
\Quad{for all}
	q \in \mbn^V.
\]
Moreover, if $Q^0$ is in equilibrium, then,
writing $\lambda_{\min} \cq \min_{v \in V} \lambda_v$,
\[
	\ex{Q^0_v}
\le
	\frac{6 n \log(2n/e)}{\beta \lambda_{\min} (s_v - \nu_v)^2}
\Quad{for all}
	v \in V.
\]
\end{introthm}

We now evaluate this under the heuristic-driven choice from \cref{cor:mix}.

\begin{introcor}[Heuristic-Driven Choice]
\label{cor:q}
Suppose that $\lambda_v = d_v / \bar d$ and $p_v \le \tfrac23 K / d_v$ for all $v \in V$, where $\bar d \cq \tfrac1n \sum_{v \in V} d_v$ is the average degree.
Let $\delta \cq \min_{v \in V} d_v$.
Suppose also that $\nu_v < s_v$ for all $v \in V$.
Let $Q \sim \QMCH_\Omega(\BM \lambda, \BM p; \BM \nu, \BM 1)$.
Then, in equilibrium,
\[
	\ex{Q^0_v}
\le
	\frac{18 \bar d n \log(2n/e)}{\delta (s_v - \nu_v)^2}
\Quad{for all}
	v \in V.
\]
\end{introcor}

A related result was proved by Jiang et al~\cite{JLNRW:parallel} for the usual hardcore model (one colour). Also, they restrict to the special case $p_v = p < 1/\Delta$, where $\Delta$ is the maximum degree.

\renewcommand{\thesubsection}{\thesection.\arabic{subsection}}

\section{Motivation and Related Work}

\subsection*{\ensuremath{-}\hspace{0.5em}Fibreoptic Routing Application}

The original motivation of this model was to create a \emph{fully decentralised} random access scheme for resource sharing in fibreoptic routing networks.
There, nodes are connected by \textit{links}, and they communicate with each other along \textit{routes}, which are sequences of links.
Multiple routes may share a subset of links; such routes \textit{interfere}.
Each link has a collection of \textit{frequencies} available.

A naive approach to routing is for the source node to send the data to the first intermediary node on the route, along with instructions of where to send on. That intermediary node processes the data and sends it onto the next node.
This continues until the data reaches its target destination.

It is possible for different frequencies to be used along the route, due to the intermediary processing.
When checking whether it is possible for a certain collection of routes to be active simultaneously, it is enough to check that no individual link is overloaded.
However, the intermediary processing adds overhead.
	If the time it takes to transmit the data along the link is larger than the processing time, then the overhead is unimportant.
However, in fibreoptic networks, data is sent along links extremely quickly, and the processing overhead becomes the performance bottleneck.

Instead of processing and resending the data at an intermediary node, an \textit{optical switch} is configured. This switch is like a prism: light coming from a single source is sent in different directions, depending on its frequency (colour). This allows a \textit{light path} to be set up along a route, removing the processing overhead; however, the \emph{same} frequency must be used throughout the entire route.

The difficulty is in choosing the frequency (colour) of the light path.
Now, it is not enough to simply check that each link is not overloaded marginally, as the colours are correlated.
In the set-up of the multicoloured hardcore model, the vertices correspond to routes, and two routes (vertices) are adjacent, forming an edge, if (and only if) they interfere---ie, share a link.
Certainly, not all routes will be able to be on simultaneously; an access scheme must be devised.

I originally learnt of this model from a talk at the Isaac Netwon Institute by \textcite{W:BT-fibreoptic} at the \href{https://www.newton.ac.uk/event/tgmw57/}{\color{black}\emph{Algorithms and Software for Quantum Computers}} event.
There, the speaker was looking to quantum computation to assist with this problem.
I, as a probabilist, took a randomised approach.

The multicoloured hardcore model has the significant benefit of \emph{decentralisation}.
All decisions made can be made by the individual vertices, without any need for synchronisation or knowledge of the state of the other routes.
A vertex can even request a light path blindly \cite{K:private-fibreoptic}:
	the path is set up if it does not conflict with any other already-active paths;
	otherwise, an error is returned to the initiator.
Moreover, the optical-switch reconfiguration is fast and easy.

The hardcore model is a popular and well-studied model for random access schemes where there is only a single frequency: \on or \off.
A toy model for this is local radio communication:
	vertices represent pairs of agents who wish to communicate;
	nearby pairs of agents cannot communicate simultaneously.
Quite separately, Glauber dynamics are used to sample proper colourings on a graph.
It seems natural to combine these two, yielding a multicoloured hardcore model which can model more complex interference situations,
	such as when multiple independent radio frequencies are available.
However, to the best of my knowledge, it has not been studied before in this context.


\subsection*{\ensuremath{-}\hspace{0.5em}Multihop Wireless Networks}

Another application of this type of random routing scheme is to \emph{multihop wireless networks}.
In cellular and wireless local area networks, wireless communication only occurs on the last link between a base station and the wireless end system.
In \emph{multihop} wireless networks, there are one or more intermediate nodes along the path; these receive and forward packets via the wireless links.
There are several benefits to the multihop approach, including extneded coverage and improved connectivity, higher transfer rates and the avoidance of wide deployment of cables.
Unfortunately, protocols, particularly those for routing, developed for fixed or cellular networks, or the Internet, are not optimal for these, more complicated, multihop wireless networks; see, eg, \cite{BN:wireless-survey}.

A highly prominent example of multihop wireless networks is in the development and deployment of 5G cellular networks \cite{TFL:multihop-5g}.
Conventional cellular networks employ well-planned deployment of tower-mounted base stations. They are undergoing a fundamental change to deployment of smaller base stations. Multihop relaying can be instrumental for implementation.
See \cite[\S 4.1]{HM:multihop-5g} for more details, from which part of this short paragraph was paraphrased.

\medskip

A multihop network with a single transmission frequency
falls precisely into the framework of the (usual) hardcore model.
Glauber dynamics is a powerful tool used to generated randomised, approximate solutions to combinatorially difficult problems.
Moreover, it often has natural decentralised implementation.
It has already been used in the past to design and analyse distributed scheduling algorithms for multi\-hop wireless networks;
see, particularly, \cite{JLNRW:parallel,BKKRSH:multihop}, from which this paragraph was paraphrased, as well as \cite{NS:dist-csma,NTS:q-csma,BKMS:multihop,JW:dist-csma,RSS:adiabatic}.

Multi\-hop wireless networks with \emph{multiple} transmission frequencies correspond precisely to our model.
To the best of our knowledge, this model has received little attention.
However, with technological and engineering advances, it could become an important extension in the future.

\subsection*{\ensuremath{-}\hspace{0.5em}(A)synchronicity}

One aspect to point out is our lack of synchronicity:
	we use \emph{continuous time}, so sites update one at a time.
In practice, engineering implementations often prefer synchronised updates.
This is the case in \cite{JLNRW:parallel}, where the (usual) hardcore model is analysed and an \emph{independent set} of vertices%
	---to recall, a set of vertices with no edges between them---%
is updated simultaneously.
It is crucial that it is an \emph{independent} set:
	that way, the changes to one vertex in the set do not affect the other vertices,
	and the updates can be done independently, in a parallel, distributed manner.

The (independent) set of vertices still needs to be chosen in each step.
In \cite{JLNRW:parallel}, the authors simply prescribe a distribution $q$ over the collection of all independent sets; no comment is made on \emph{how} to draw from this set.
In principle, this distribution is very complicated, and perhaps even needs approximating---eg, via Glauber dynamics for the (usual) hardcore model.

The path coupling technique that we use, and is used in \cite{JLNRW:parallel}, is robust to parallel updates, provided one update does not affect the other updates,
	such as is the case for updating an independent set of vertices.
If $N$ is the expected cardinality of the independent set chosen%
	---ie, $N \cq \ex[S \sim q]{\abs S}$---%
then the mixing bound behaves as if time is sped up by a factor $N$.
We consider single-site, continuous-time updates for simplicity; our analysis extends to the parallel set-up, too.

\subsection*{\ensuremath{-}\hspace{0.5em}Spin Systems in Statistical Mechanics}

Spin systems are widely studied in statistical mechanics,
crossing combinatorics, probability and physics:
	these involve a graph $G = (V, E)$ and a discrete set $\mck$ of \textit{spins};
	each vertex $v \in V$ is assigned a spin $k \in \mck$.
Adjacent vertices interact with each other.
A whole zoo of examples of discrete spin systems is discussed extensively in the very recent paper by \textcite{PS:spin-systems}.

\begin{itemize}
	\item 
	In proper colourings, $\mck = \bra{1, ..., K}$ and the constraint is hard:
		adjacent vertices \emph{must not} have the same colour.
	The hardcore model is similar with $\mck = \bra{0,1}$.
	
	\item 
	In the Ising model, $\mck = \bra{\pm1}$ and the constraint is soft:
		vertices prefer to be aligned with their neighbours,
		with strength controlled by the \textit{inverse temperature} $\beta \ge 0$.
\end{itemize}
The multicoloured hardcore model is discussed, under the name \textit{anti-Widom--Rowlinson}, in \cite[\S 3.2.2]{PS:spin-systems}.
It was originally introduced by \textcite{RL:mch} in the context of lattice~gases.

The results of \cite{RL:mch,PS:spin-systems} are specialised to the lattice $\mbz^d$.
The latter is most interested in the case where the dimension $d$ is much larger than the number $K$ of colours.
The motivating example for this paper is the fibreoptic routing, for which the lattice $\mbz^d$---particularly in high dimensions---is not an appropriate model.
Our results appear to be the first on general graphs.

\subsection*{\ensuremath{-}\hspace{0.5em}Notation}

We briefly recall some notation which is used throughout the paper.

\begin{itemize}
\item 
The underlying graph is $G = (V, E)$.
Let $n \cq \abs V$ denote its number of vertices, and write $u \sim v$ if $\bra{u,v} \in E$.
The degree of vertex $v \in V$ is $d_v \cq \sum_{u \in V} \one{ u \sim v } = \abs{ \bra{ u \in V \mid u \sim v } }$.

\item 
There are $K \in \mbn$ colours, and we abbreviate $[K] \cq \bra{1, ..., K}$ and $[K]_0 \cq \bra{0, 1, ..., K}$.

\item 
The update rates and probabilities are $\BM \lambda \in (0, \infty)^V$ and $\BM p \in [0,1]^V$, respectively.

\item 
The state space is $\Omega \cq \bra{ \omega \in [K]_0^V \mid \omega \text{ is proper} }$, where $\omega \in [K]_0^V$ is \textit{proper} if
\[
	\omega_u \ne \omega_v
\Quad{whenever}
	\bra{u,v} \in E
\Quad{and}
	\omega_u + \omega_v > 0.
\]

\item 
The multicoloured hardcore model is denoted $\MCH_\Omega(\BM \lambda, \BM p)$, and its equilibrium distribution~$\pi$.

\item 
For the queue $\QMCH_\Omega(\BM \lambda, \BM p; \BM \nu, \BM 1)$, the
arrival rates are $\BM \nu \in (0, \infty)^V$ and services rates uniformly (over $v \in V$) are $1$.
and the \textit{equilibrium service rates} are $\BM s = (s_v)_{v \in V}$, where
\[
	s_v
\cq
	\sumt{\omega \in \Omega : \omega_v \ne 0}
	\pi(\omega)
\Quad{for}
	v \in V.
\]
\end{itemize}

\renewcommand{\thesubsection}{\Alph{subsection}}
\section{Proofs of Main Theorems}

\subsection{Mixing}

In this section, we use the classical path coupling argument of \textcite{BD:path-coupling} to upper bound the mixing time.
Throughout, $X, Y \sim \MCH_\Omega(\BM \lambda, \BM p)$, under the `natural' coupling:
\begin{itemize}
	\item 
	the update clocks of each vertex are coupled, so the same vertex is chosen at the same time;
	
	\item 
	the subsequent coin toss and colour selection are also coupled.
\end{itemize}
This coupling is clearly \textit{coalescent}:
\[
	X^t = Y^t
\Quad{implies}
	X^s = Y^s
\Quad{for all}
	s \ge t.
\]

%
%
%


\begin{Proof}[Proof of \cref{thm:mix}]
We use path coupling, so must define a path space.
We say that $x,y \in [K]_0^V$ are \emph{adjacent} if there is a unique $v \in V$ such that $x_v \ne y_v$ \emph{and} $0 \in \bra{x_v, y_v}$.
In other words, our path space is generated by
	activating an inactive vertex
or
	deactivating an active vertex;
changing the colour of an already active vertex \emph{is not} permitted.
This space is connected:
	let $d(x,y)$ denote the distance between two configurations $x,y \in [K]_0^V$;
	then, $\one{x \ne y'} \le d(x, y) \le 2n$ for all $x,y \in [K]_0^V$,
	going via the empty configuration $(0,...,0) \in \Omega$.

For $v \in V$ and $x \in [K]_0^V$,
denote the \textit{available colours at $v$ in $x$} by
\[
	\mca_v(x)
\cq
	[K] \setminus \cup_{u \in V : \bra{u,v} \in E} \bra{x_u}
=
	\bra{ k \in [K] \mid x_u \ne k \: \forall \, u \sim v }.
\]

Suppose that $(X^0, Y^0) = (x,y) \in \Omega^2$ with $d(x,y) = 1$; say, $0 = x_v \ne y_v$.
Consider the first step of the process from these states.
Suppose that vertex $u \in V$ updates,
which happens at rate $\lambda_u$.

\begin{itemize}
	\item 
	Suppose that $u \not\sim v$.
	Then, $\mca_x(u) = \mca_y(u)$, since $x_w = y_w$ for all $w \sim u$.
	Hence, we can perform the same update in both $X$ and $Y$.
	The relative distance is unchanged, unless $u = v$, in which case the two coalesce.
	
	\item 
	Suppose that $u \sim v$; in particular, $u \ne v$.
	We may not have $\mca_u(x) = \mca_u(y)$, but we do have
	\[
		\mca_u(x) \cup \bra{x_v}
	=
		\mca_u(y) \cup \bra{y_v};
	\Quad{but,}
		x_v = 0 \notin \mca_u(x) \subseteq [K].
	\]
	Hence, $\abs{ \mca_u(x) \mathrel{\triangle} \mca_u(y) } \le 1$.
	So, the probability that a proposed colour is valid for one and not the other is at most $1/K$.
	If this is the case, then the relative distance increases by $1$; otherwise, it remains unchanged.
	The probability that \emph{some} colour is proposed is $p_u$.
\end{itemize}
It is in this last step that the assumption $0 \in \bra{x_v, y_v}$ is used:
	without it, the symmetric difference could be of size $2$, giving a probability $2/K$.
Summing over $u \in V$, weighted by $\lambda_u$, the relative distance increases by $1$ at rate at most $\tfrac1K \sum_{u : u \sim v} p_u \lambda_u$ and decreases by $1$ at rate $\lambda_v$.
Hence,
\[
	\left.
	\tfrac d{dt}
	\ex[x,y]{ d(X^t, Y^t) }
	\right|_{t=0}
\le
	\lambda_v
	\rbb{ \tfrac1K \sumt{u : u \sim v} p_u \lambda_u / \lambda_v - 1 }
\le
	- \beta \lambda_v,
\]
with the last inequality using the (main) assumption of the theorem.
This can be extended to general $x,y \in [K]_0^V$%
	---ie, not requiring $d(x,y) = 1$---%
by looking at contraction along geodesics, in the usual manner for path coupling.
Hence, recalling that $\lambda_{\min} = \min_v \lambda_v$,
\[
	\max_{x,y \in [K]_0^V}
	\left.
	\tfrac d{dt}
	\exb[x,y]{ d(X^t, Y^t) }
	\right|_{t=0}
\le
	- \beta \lambda_{\min}.
\]
By the Gr\"onwall inequality, integrating this and using $\one{x \ne y} \le d(x,y) \le 2n$, we obtain
\[
	\max_{x,y \in [K]_0^V}
	\pr[x,y]{ X^t \ne Y^t }
\le
	\max_{x,y \in [K]_0^V}
	\ex[x,y]{ d(X^t, Y^t) }
\le
	2 n e^{-\beta t}.
\]
Finally, the coupling representation of total-variation distance implies that
\[
	\max_{x,y \in \Omega}
	\tv{ \pr[x]{ X^t \in \cdot } - \pr[y]{ Y^t \in \cdot } }
\le
	\min\brb{ 2 n e^{- \beta \lambda_{\min} t}, \: 1 }.
\qedhere
\]
\end{Proof}

\begin{rmkt*}
If preferred, instead of using a continuous-time version of path coupling, discretise time:
	let $\tilde X^\ell \cq X^{\delta \ell}$ and $\tilde Y^\ell \cq Y^{\delta \ell}$,
	where $\delta$ is some very small real number.
Then,
\[
	\ex[x,y]{ d(\tilde X^1, \tilde Y^1) }
\le
	\rbb{ 1 - \beta \lambda_{\min} \delta + \oh \delta }
	d(x,y)
\quad
	\text{uniformly},
\]
using the fact that the diameter is finite to obtain a uniform $\oh \delta$ term.
Then, path coupling gives
\[
	\ex[x,y]{ d(\tilde X^\ell, Y^\ell) }
\le
	2n
	\rbb{ 1 - \beta \lambda_{\min} \delta + \oh \delta }{}^\ell
\le
	2n
	e^{-\beta \lambda_{\min} \delta \ell + \oh{\delta \ell}} n.
\]
Given $t \ge 0$, let $\ell \cq \floor{t/\delta} \ge t/\delta - 1$.
Then,
\[
	\ex[x,y]{ d(X^t, Y^t) }
\le
	\ex[x,y]{ d(\tilde X^\ell, \tilde Y^\ell) }
\le
	2n
	e^{-\beta t + \oh1}.	
\]
Finally, taking $\delta \downarrow 0$, we deduce the same bound as before.
\end{rmkt*}

We close this section with a discussion of the equilibrium service rates.
Here, we assume that
\[
	p_v
\le
	\tfrac13 K / \tilde d_v
\Quad{where}
	\tilde d_v
\cq
	\max\bra{ d_u \mid u \sim v \text{ or } u = v }
\Quad{for}
	v \in V.
\]

\begin{Proof}[Proof of \cref{prop:rates}]
The most important quantity to estimate is the proportion of colours available at a vertex.
This allows estimation of the probability an attempted colouring is successful.

Let $v \in V$.
Clearly, in equilibrium, each neighbour $u$ of $v$ is active with probability at most $p_u \le \tfrac13 K / \tilde d_u$; in particular, $s_v \le p_v$.
Hence, if $N_v$ is the number of colours available at $v$, then
\[
	N_v
\lesssim
	\Bin(d_v, \tfrac13 K / d_v)
\quad
	\text{in equilibrium}.
\]
It can be shown that
\(
	\pr{ \Bin(d, \tfrac13 k / d) \ge \tfrac12 k }
\le
	\tfrac13
\)
whenever $k \le d$.
This implies that
\[
	\pr{ N_v \ge \tfrac12 K }
\le
	\tfrac13.
\]
Hence, upon refreshing, at least $\tfrac12$ of the colours are available with probability at least $\tfrac23$.
In particular, the probability that the proposed colour is accepted is at least $\tfrac13$.
Thus,
\(
	s_v
\ge
	\tfrac13
	p_v.
\)
\end{Proof}

We discuss briefly extensions of this proof, including heuristics for an upper bound on $s_v$.

\begin{rmkt*}
If we require $p_v \le (1 - \delta) K / \tilde d_v$, then the above argument says that at least a proportion $\delta$ of the colours are free in expectation.
If $K$ (and $\tilde d_v$) are large, then the Binomial concentrates around its expectation.
There is then a probability $\delta$ that a uniformly proposed colour is available.

We can extend this, heuristically at least.
If $u,u' \sim v$, then the colours at $u$ and $u'$ should be approximately independent if $K$ is large and the graph has few triangles. If $k_1, ..., k_K \sim^\iid \Unif([K])$, then $\tfrac1K \abs{ \bra{ k_1, ..., k_K } } \approx 1/e$, suggesting that, in fact, a proportion $1/e$ are available after $K$ choices.

We can also try to iterate this argument.
Instead of upper bounding the expected number of colours taken by $\sum_{u : u \sim v} p_u$, we can upper bound by $\sum_{u : u \sim v} s_u$.
Suppose that $s_v$ does not vary much over the vertices: $s_v \approx \bar s \cq \tfrac1n \sum_u s_u$, the average of $s$; see, eg, Figure~\ref{fig:rr} later.
Also, assume graph regularity: $d_v = d$, and $p_v = p$, for all $v$.
Then, $\sum_{u : u \sim v} s_u \approx d \bar s$.
This imposes
\[
	\bar s
\le
	p(1 - d \bar s / K);
\Quad{ie,}
	\bar s
\le
	p / (1 + p d / K).
\]
Including the factor $1/e$ from the previous heuristic improves this to
\(
	\bar s
\approx
	p / \rbr{ 1 + e^{-1} p d / K }.
\)
\end{rmkt*}

\subsection{Queues}

In this section, we investigate the \textit{stability} of the queueing network:
	ie, its positive recurrence (or lack thereof) and expected queue length in equilibrium.
The end goal is \cref{thm:q}.
Similar properties for a related model are established in \cite[\S V]{JLNRW:parallel}, using the standard \textit{Lyapunov function}
\[
	L^t
\cq
	\sumt{v \in V}
	(Q^t_v)^2
\Quad{for}
	t \ge 0
\Quad{where}
	Q = (Q^t)_{t\ge0}
\sim
	\QMCH_\Omega(\BM \lambda, \BM p; \BM \nu, \BM 1).
\]
There, the model is slightly simpler, with unit service times, rather than the usual Exponentials.
Moreover, they require $p_v = p \le 1/\Delta$ for all $v \in V$, where $\Delta \cq \max_v d_v$ is the maximum degree of the graph $G = (V, E)$, and treat $\Delta$ as a constant, which is absorbed into a final, unquantified constant.
For a sequence $(G_n)_{n\in\mbn}$ of graphs, this implicitly assumes bounded degrees: $\sup_{n\in\mbn} \Delta_n < \infty$.
We allow much greater generality, both in the graph and the choice of $\BM p = (p_v)_{v \in V}$.

We denote by $\tau$ the first time the queue is empty:
\[
	\tau
\cq
	\inf\brb{ t \ge 0 \mid Q^t = 0, \: \cup_{s \in [0,t]} \, Q^s \ne \bra{0} }.
\]
Positive recurrence is equivalent to having $\ex[q]{\tau} < \infty$ for some, and hence all, $q \ne 0$.

%
%
%
%
%

\begin{Proof}[Proof of \cref{thm:q}]
We define an appropriate Lyapunov function $L$, and establish negative drift:
\[
	L^t
\cq
	\tfrac12
	\sumt{v \in V}
	(Q^t_v)^2
\Quad{for}
	t \ge 0.
\label{eq:q:L-def}
\nt
\]

We fix some notation and conventions.
By the memoryless property of the service times, we may assume that the vertices are \emph{always} providing service, but that a service attempt is rejected if the vertex is inactive (ie, has colour $0$) at the time of the attempt.
Then, the arrivals and attempted services form Poisson processes, independent of each other and of the underlying \MCH process.

Fix $v \in V$ and $t,T \ge 0$.
Write $\hat S_v[T, T+t)$ for the number of attempted services by vertex $v$ between times $T$ and $T+t$, and write $\hat s_v \cq \hat s_v[T, T+t) \cq \tfrac1t \hat S_v[T, T+t)$ for the average (attempted) service rate in this interval.
Similarly, write $\hat A_v[T, T+t)$ and $\hat a_v \cq \hat a_v[T, T+t) \cq \tfrac1t \hat A_v[T, T+t)$ for the number of arrivals and average service rate, respectively, between $T$ and $T+t$.

Using these definitions, we have the following simple inequality:
\[
	Q^{T+t}_v
\le
	\sbb{ Q^t_v - \hat S_v[T, T+t) }_+
+	\hat A_v[T, T+t)
=
	\sbr{ Q^t_v - t \hat s_v }_+
+	t \hat a_v,
\]
where $[\alpha]_+ \cq \max\bra{\alpha, 0}$ for $\alpha \in \mbr$.
Hence,
using $[Q^t_v - t \hat s_v]_+ \le Q^t_v$,
\[
\begin{aligned}
	(Q^{T+t}_v)^2
&
\le
	\rbr{ Q^T_v - t \hat s_v }^2
+	2 t \sbr{ Q^T_v - t \hat s_v }_+ \hat a_v
+	t^2 \hat a_v^2
\\&
\le
	(Q^T_v)^2
+	2 t Q^T_v (\hat a_v - \hat s_v)
+	t^2 (\hat a_v^2 + \hat s_v)^2.
\end{aligned}
\label{eq:q:Q2-bound}
\nt
\]
Plugging this into the definition \eqref{eq:q:L-def} of $L$ bounds its random increment:
\[
	L^{T+t} - L^T
\le
	t
	\sumt{v \in V}
	Q^T_v \rbr{ \hat a_v - \hat s_v }
+	\tfrac12 t^2
	\sumt{v \in V}
	\rbr{ \hat a_v^2 + \hat s_v^2 }.
\label{eq:q:L-incr}
\nt
\]

Now, if $\hat \tau_v$ is the proportion of time during $[T, T+t)$ that vertex $v$ is active, then
\[
	t \hat a_v
=
	\hat A_v[T, T+t)
\sim
	\Pois(t \nu_v)
\Quad{and}
	t \hat s_v
=
	\hat S_v[T, T+t)
\sim
	\Pois(t \hat \tau_v).
\]
To emphasise, the implicit Poisson variables are independent of the \MCH process.
Recall that
\[
	\text{if}
\quad
	P \sim \Pois(m),
\Quad{then}
	\ex{P} = m
\Quad{and}
	\ex{P^2} = m + m^2.
\]
Now, $\nu_v < s_v$, by assumption, and $s_v \le p_v \le 1$; also, $\hat \tau_v \le 1$.
Hence,
\[
	\ex{\hat a_v}
=
	\nu_v,
\quad
	\ex{\hat a_v^2}
\le
	2,
\quad
	\ex{\hat s_v}
\le
	1
\Quad{and}
	\ex{\hat s_v^2}
\le
	2.
\]
Plugging these into \eqref{eq:q:L-incr} bounds the (expected) drift:
\[
	\ex{ L^{T+t} - L^T \mid (X^T, Q^T) }
\le
	t
	\sumt{v \in V}
	Q^T_v \rbr{ \nu_v - \ex{\hat s_v \mid X^T} }
+	\tfrac32 n t^2;
\label{eq:q:L-drift1}
\nt
\]
the (attempted) service rate $\hat s_v[T, T+t)$ depends only on $X^T$, not $Q^T$.

It remains to handle $\ex{ \hat s_v \mid X^T }$.
The attempted services are a thinned Poisson process.
So,
\[
	\ex{ \hat s_v \mid X^T }
=
	\ex{ \hat \tau_v \mid X^T }
\Quad{and}
	\tau_v
=
	\tfrac1t
	\sumt[T+t]{T}
	\one{ X^s_v \ne 0 }
	ds.
\]
So, if we write $\mu_{x,s}$ for the law of $X^s$ given $X^0 = x$, then
\[
	\ex{ \hat s_v \mid X^T }
=
	\tfrac1t
	\intt[T+t]{T}
	\pr{ X^s \ne 0 \mid X^T }
	ds
	\tfrac1t
	\intt[t]{0}
	\mu_{X^T,s}\rbr{ \bra{ \omega \in \Omega \mid \omega_v \ne 0 } }
	ds.
\]
This is very similar to the equilibrium (attempted) service rate
\[
	s_v
=
	\sumt{\omega \in \Omega : \omega_v \ne 0}
	\pi(\omega)
=
	\pi\rbr{ \bra{ \omega \in \Omega \mid \omega_v \ne 0 } };
\]
in fact, by the ergodic theorem, $\hat s_v[T, T+t) \to s_v$ as $t \to \infty$.
Quantitatively,
\[
	\abs{ \ex{ \hat s_v \mid X^T } - s_v }
&
=
	\absb{ 
		\tfrac1t
		\intt[t]{0}
		\mu_{X^T,s}\rbr{ \bra{ \omega \in \Omega \mid \omega_v \ne 0 } }
		ds
-		\pi\rbr{ \bra{ \omega \in \Omega \mid \omega_v \ne 0 } }
	}
\\&
\le
	\tfrac1t
	\intt[t]{0}
	\abs{ 
		\mu_{X^T,s}\rbr{ \bra{ \omega \in \Omega \mid \omega_v \ne 0 } }
	-	\pi\rbr{ \bra{ \omega \in \Omega \mid \omega_v \ne 0 } }
	}
\\&
\le
	\tfrac1t
	\intt[t]{0}
	\tv{ \mu_{X^T,s} - \pi }
	ds.
\]

It is here that we apply the mixing result, \cref{thm:mix}:
for any $x \in \Omega$ and $s \ge 0$,
\[
	\tv{ \mu_{x,s} - \pi }
\le
	\min\bra{ 2n e^{-\beta \lambda_{\min} s}, 1 };
\]
note that the first hypothesis of \cref{thm:q} is precisely that required to apply \cref{thm:mix}.
Then,
\[
	\intt[t]{0}
	\tv{ \mu_{X^T,s} - \pi }
	ds
&
\le
	t_0
+	n
	\intt[t \vee t_0]{t_0}
	e^{-\beta \lambda_{\min} s}
	ds
\\&
\le
	t_0 + (\beta \lambda_{\min})^{-1}
\eqqcolon
	t_1
\Quad{where}
	t_0
\cq
	(\beta \lambda_{\min})^{-1} \log(2n).
\]
In particular, this is independent of $t$, so vanishes once divided by $t$ and $t \to \infty$:
\[
	\abs{ \ex{ \hat s_v \mid X^T } - s_v }
\le
	t_1 / t
\to
	0
\Quad{as}
	t \to \infty.
\]

Recall that we want to plug this bound into \eqref{eq:q:L-drift1}.
Let
\(
	\eps_v
\cq
	\tfrac12 \rbr{ s_v - \nu_v }
\)
and $t_v \cq t_1 / \eps_v$.
Then,
\[
	\abs{ \ex{ \hat s_v \mid X^T } - s_v }
\le
	\eps_v
\Quad{whenever}
	t \ge t_v
\Quad{for all}
	v \in V.
\]
Set $t_\star \cq \max_v t_v$, so $t_\star \ge t_v$
Plugging this into \eqref{eq:q:L-drift1},
\[
	\ex{ L^{T+t} - L^T \mid (X^T, Q^T) }
&
\le
	- t
	\sumt{v \in V}
	Q^T_v (s_v - \nu_v - \eps_v)
+	\tfrac32 n t^2
\\&
\le
	- \tfrac12 t
	\sumt{v \in V}
	Q^T_v (s_v - \nu_v)
+	\tfrac32 n t^2
\Quad{whenever}
	t \ge t_\star.
\label{eq:q:L-drift2}
\nt
\]
The last expression is negative for large enough $\norm{ Q^T }$.
This establishes negative drift of $L$.
Hence, by the Foster--Lyapunov criterion (see, eg, \cite[Proposition~D.1]{KY:stoc-net-book}), $(Q^t)_{t\ge0}$ is positive recurrent.

It remains to control the expected queue length in equilibrium.
To do this, we start in equilibrium and take the expectation of the increment $(Q^{t_v}_v)^2 - (Q^0_v)^2$.
By stationarity and \eqref{eq:q:Q2-bound},
\[
	0
=
	\ex{ (Q^{t_v}_v)^2 - (Q^0_v)^2 }
\le
	- t_v
	\ex{Q^0_v}
	\rbr{ s_v - \nu_v - \eps_v }
+	\tfrac32 n t_v^2,
\]
using the same manipulations as before.
Rearranging,
\[
	\ex{Q^0_v}
\le
	\tfrac32 n t_v / (s_v - \nu_v - \eps_v)
\le
	6 n t_1 / (s_v - \nu_v)^2.
\]
Finally, $t_1 = (\beta \lambda_{\min})^{-1} (\log(2n) + 1) = (\beta \lambda_{\min})^{-1} \log(2n/e)$.
\end{Proof}

\section{Simulations: Queue Lengths and Equilibrium Service Rate}

We close the paper with a short discussion of some simulations.
Specifically, we investigate the queue lengths and the proportion of time that a vertex is active as a rolling average---namely,
\[
\textstyle
	\hat Q^t_v
\cq
	\tfrac1t \sum_{s=0}^{t-1} Q^s_v
\Quad{and}
	\hat s^t_v
\cq
	\tfrac1t \sum_{s=0}^{t-1} \one{X^s_v \ne 0}
\Quad{for}
	t \ge 0.
\]
Then, $\hat Q^t_v \to \ex[\pi]{Q^0_v}$ and $\hat s^t_v \to s_v$ as $t \to \infty$,
the expected equilibrium queue length and service~rate.

Our choice of parameters is driven by the same heuristics as for \cref{cor:mix,cor:q}:
\[
	\lambda_v \cq d_v / \bar d,
\quad
	p_v \cq \min\bra{\tfrac45 e K / d_v, \tfrac12},
\quad
	\nu_v \cq \tfrac13 p_v
\Quad{and}
	\mu_v \cq 1
\Quad{for}
	v \in V.
\]
Notice the prefactor in $p_v$: it is $\tfrac45 e > 2$, rather than $\tfrac13$ or $\tfrac23$.
This is to emphasise the fact that we really can take $p_v$ close to $e K/d_v$, yet still get high, and stable, service rates $s_v$.


\begin{figure}[t]
\centering
\begin{subfigure}{.475\textwidth}
	\centering
	\includegraphics[width=0.98\linewidth]{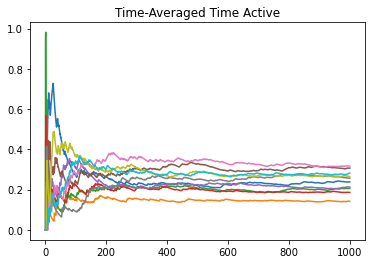}
	\label{er:active}
\end{subfigure}%
\hfill
\begin{subfigure}{.475\textwidth}
	\centering
	\includegraphics[width=0.98\linewidth]{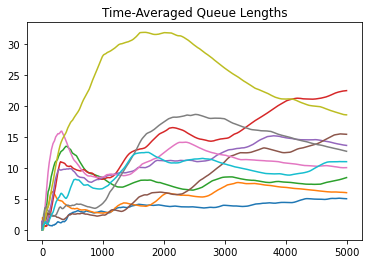}
	\label{er:queues}
\end{subfigure}%
\caption{The underlying graph is Erd\H{o}s--R\'enyi with $n = 500$ vertices and edge probability~$40/n$}
\label{fig:er}
\vspace{-\medskipamount}
\end{figure}

\begin{figure}[t]
\centering
\begin{subfigure}{.475\textwidth}
	\centering
	\includegraphics[width=0.98\linewidth]{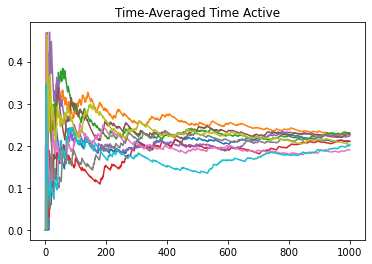}
	\label{cm:active}
\end{subfigure}%
\hfill
\begin{subfigure}{.475\textwidth}
	\centering
	\includegraphics[width=0.98\linewidth]{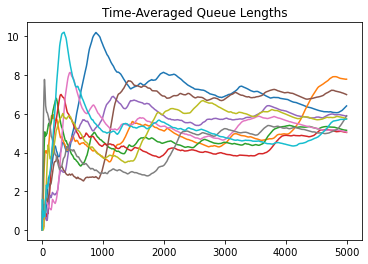}
	\label{cm:queues}
\end{subfigure}%
\caption{The underlying graph is drawn uniformly over $40$-regular graphs on $n = 500$ vertices}
\label{fig:rr}
\vspace{-\medskipamount}
\end{figure}

The plots in
Figure~\ref{fig:er} show the time-averaged queue lengths and service rates when the underlying graph is an Erd\H{o}s--R\'enyi graph.
Those in Figure~\ref{fig:rr} show the same for a random regular graph.
The average degree is $40$ and $K = 10$ colours are used;
so, almost all vertices satisfy $p_v = \tfrac45 e K / d_v \approx 0.5$.
A collection of $10$ vertices with typical degrees to be displayed are chosen randomly.
Time is scaled so that the average vertex update-rate is $1$%
	---ie, scaled by $\tfrac1n \sum_v (\lambda_v + \nu_v + 1)$.

We see that the empirical service rates settle down really quite quickly, and appear to be remain stable.
Moreover, the values $s_v$ to which they converge appear to be on the same order as the proposal probabilities $p_v$.
This suggests that many proposals are accepted, but not \emph{too many}: if $s_v \approx p_v$, then this suggests that a higher proposal probability $p_v$ could have been used.
In particular, we found that the normalised difference
\(
	\abs{ s_v - p_v } / p_v
\)
averaged around 60\% (over $v \in V$).

The queue lengths, on the other hand, fluctuate a more.
They are a bit more stable in the random regular graph (Figure~\ref{fig:rr}) compared with the Erd\H{o}s--R\'enyi graph (Figure~\ref{fig:er}).
This is perhaps due to the inhomogeneities in the latter.
It is not even completely clear what they are converging~to.


We suggest that this is likely caused by the inhomogeneities in the graph along with the fact that we take $\nu_v = p_v / 3 \approx 0.33 p_v$, which is pretty close to $s_v \approx 0.4 p_v$.
Indeed, the same calculations (not shown) with $\nu_v = 0.2 p_v$ result in much more stable queue lengths.

The primary objective is to get as large an equilibrium service rate $s_v$ as possible, or at least a large average $\bar s = \tfrac1n \sum_v s_v$.
Since the 60\% above is still quite a large rejection rate, we also tested a slightly smaller value of $p_v$: namely, we used $p_v = \tfrac23 e K / d_v \approx 0.45$.
However, we found that $\bar s$ was about 10\% smaller for these parameters, for both random graph models.

%

A random $d$-regular graph locally looks like a $d$-regular tree, so it is not reasonable to expect better than $e K / d = p_c(\Delta, K)$, the earlier critical threshold.
Similarly, a sparse Erd\H{o}s--R\'enyi graph locally looks like a Bienaymé--Galton--Watson tree with $\Pois(\bar d)$ degrees, with mean $\bar d$.

\vfill

\noindent%
\textbf{\sffamily Acknowledgements.}
The author would like to thank Frank Kelly for multiple detailed discussions on this topic, as well as reading through the paper.
Additionally, thanks go to Perla Sousi and Luca Zanetti, who also read the paper.
Their feedback has been invaluable in preparing this manuscript.

\renewcommand{\bibfont}{\sffamily\small}
\printbibliography[heading=bibintoc]

\end{document}